\documentclass[11pt]{amsart}
\usepackage{amssymb,amsmath,mathtools}
\usepackage[margin=1.15in]{geometry}
\usepackage[hidelinks]{hyperref}
\newtheorem{theorem}{Theorem}
\newtheorem{lemma}[theorem]{Lemma}
\newtheorem{corollary}[theorem]{Corollary}
\newtheorem{proposition}[theorem]{Proposition}
\theoremstyle{definition}
\newtheorem{conjecture}[theorem]{Conjecture}
\theoremstyle{remark}
\newtheorem{remark}[theorem]{Remark}

\newcommand{\cL}{\mathcal{L}}
\newcommand{\cW}{\mathcal{W}}
\newcommand{\A}{\mathcal{A}}

\title{Counting words of polynomially bounded complexity}
\author{Nick Stucky-Mack}
\address{Symbolic Mirror}
\email{nick@symbolicmirror.info}
\date{\today}

\begin{document}

\begin{abstract}
Let $\cL_n(f)$ be the set of length-$n$ factors of infinite words whose complexity function is bounded by $f$. For $f(n)=n^\alpha$ with $\alpha>1$, Mauduit and Moreira showed that $\log|\cL_n(f)|$ lies between $c\,n^{1-1/\alpha}/\log n$ and $c'(n\log n)^{(\alpha+1)/(\alpha+2)}$. We prove $\log|\cL_n(Cn^\alpha)|=\Theta(n^{1-1/\alpha}\log n)$ for $C\ge2$. The upper bound refines a Lempel--Ziv lemma of Raskhodnikova, Ron, Rubinfeld and Smith. For linear $f$, that lemma and a construction of Kociumaka, Navarro and Prezza give $\log|\cL_n(f)|=\Theta(\log^2n)$ when $f(n)\ge2n$, which answers both questions of Mauduit and Moreira on linear complexity. In both cases we determine the critical gauge for the generalized Hausdorff measure of the set of reals whose base-$q$ expansion has complexity bounded by $f$, up to a constant in the exponent; for $2n\le f(n)\le An$ this answers, for base-$q$ expansions, a further question of Mauduit and Moreira. We also note the consequence for correlations of random signs and of the M\"obius function with words of bounded complexity.
\end{abstract}

\maketitle

\section{Setting and results}

Let $\A$ be a finite alphabet with $q\ge2$ letters. For a word $u$ (finite or infinite) let $d_\ell(u)$ be the number of distinct factors of $u$ of length $\ell$; for infinite $w$, $p_w=d_\cdot(w)$ is the complexity function. Following Mauduit and Moreira~\cite{MM10}, for $f:\mathbb N\to\mathbb R^+$ let
\[
\cW(f)=\{w\in\A^{\mathbb N}: p_w(\ell)\le f(\ell)\ \forall\ell\},\qquad \cL_n(f)=\bigcup_{w\in\cW(f)}\{\text{factors of }w\text{ of length }n\},
\]
and $W_n(f)=\{u\in\A^n: d_\ell(u)\le f(\ell)\ \forall\ell\le n\}$. Since $d_\ell(u)\le p_w(\ell)$ for every factor $u$ of $w$,
\begin{equation}\label{eq:incl}
\cL_n(f)\subseteq W_n(f).
\end{equation}

For $f(n)=n^\alpha$ ($\alpha>1$, modified for small $n$), \cite[Examples~4.14]{MM10} gives
\[
\exp\big(c\,n^{(\alpha-1)/\alpha}/\log n\big)<|\cL_n(f)|<\exp\big(c'(n\log n)^{(\alpha+1)/(\alpha+2)}\big).
\]
These bounds are restated in \cite[Theorem~A]{MM12}, where they are used to estimate generalized Hausdorff dimensions.

\begin{theorem}\label{thm:main}
Let $\alpha>1$.
\begin{enumerate}
\item[(a)] If $f(\ell)\le C\ell^\alpha$ for all $\ell$, then $\log|W_n(f)|\le K\,n^{1-1/\alpha}\log n$ for $n\ge2$, where $K=K(\alpha,C,q)$.
\item[(b)] If $f(\ell)=C\ell^\alpha$ with $C\ge2$, then $\log|\cL_n(f)|\ge\big(\tfrac1{2\alpha}-o(1)\big)\,n^{1-1/\alpha}\log n$.
\end{enumerate}
Hence $\log|\cL_n(Cn^\alpha)|=\Theta(n^{1-1/\alpha}\log n)$ for $C\ge2$.
\end{theorem}

For linear bounds the corresponding statement follows from known results; we record it because it settles \cite[\S4.4]{MM10}.

\begin{corollary}\label{cor:linear}
\begin{enumerate}
\item[(a)] If $f(\ell)\le A\ell$ for all $\ell$, then $\log|W_n(f)|\le K_{A,q}\log^2n$ \textup{\cite{RRRS,KNP}}.
\item[(b)] If $f(\ell)\ge2\ell$ for all $\ell$, then $\log|\cL_n(f)|\ge c\log^2n$ \textup{\cite[\S IV]{KNP}}.
\item[(c)] If $f(\ell)\ge\ell+1$ for all $\ell$, then $|\cL_n(f)|\ge(1+o(1))n^3/\pi^2$ \textup{\cite{Lip82,Mig91}}.
\end{enumerate}
\end{corollary}

Section~4 answers the questions of \cite[\S4.4]{MM10}. Proposition~\ref{prop:hausdorff} uses these counts to determine, up to a constant in the exponent, the critical gauge for the generalized Hausdorff measure of the reals whose base-$q$ expansion lies in $\cW(f)$. Section~6 gives the consequence for correlations with random signs and with the M\"obius function.

\section{Lempel--Ziv phrases and factor counts}

The \emph{LZ77 parsing} of $u=u_1\cdots u_n$ splits $u$ greedily into phrases. If a phrase starts at $s$, let $L\ge0$ be maximal such that $u[s,s+L)$ has an occurrence starting at some $t<s$ (overlaps allowed). The phrase is then $u[s,s+L]$: a copied block of length $L$ followed by one literal letter (the last phrase may lack the literal). Let $z(u)$ be the number of phrases. The triples $(t,L,u_{s+L})$ determine $u$, so
\begin{equation}\label{eq:encode}
\#\{u\in\A^n:z(u)\le Z\}\le2\big((q+1)(n+1)^2\big)^Z.
\end{equation}

The following is \cite[Lemma~3]{RRRS} (see also \cite[Lemma~II.2]{KNP}) in dyadic form. We include the proof because Theorem~\ref{thm:main}(a) needs the per-scale bound.

\begin{lemma}\label{lem:phrases}
For $u\in\A^n$ and $r\ge0$, the number of phrases with copied length $L\in[2^r,2^{r+1})$ is at most $d_{2^{r+2}}(u)/2^r+6$.
\end{lemma}

\begin{proof}
Let such a phrase start at $s$, with literal at $s+L$. By maximality of $L$, $v=u[s,s+L]$ has no occurrence starting before $s$. For $0\le j<2^r$ the window $W_j=u[s-j,s-j+2^{r+2})$ contains $v$ at offset $j$, because $j+L+1<3\cdot2^r$. If $W_j$ occurred at $t'<s-j$, then $v$ would occur at $t'+j<s$. So $s-j$ is the first occurrence of $W_j$, and distinct such positions give distinct factors of length $2^{r+2}$.

Starts of these phrases are more than $2^r$ apart, so the sets $\{s-j:0\le j<2^r\}$ are disjoint. Three kinds of phrase are excluded: at most one with $s\le2^r$ (where $s-j<1$ is possible); at most four with $s>n-2^{r+2}+1$ (where $W_j$ may leave $u$); and the final phrase. Every other phrase contributes $2^r$ first-occurrence positions.
\end{proof}

Phrases with $L=0$ introduce a new letter, so there are at most $q$ of them. Phrases with $L\ge2^R$ number at most $n/2^R$. Hence, for every $R\ge0$,
\begin{equation}\label{eq:z}
z(u)\le q+\frac n{2^R}+\sum_{r=0}^{R-1}\Big(\frac{d_{2^{r+2}}(u)}{2^r}+6\Big).
\end{equation}

\section{Proofs}

\begin{proof}[Proof of Theorem~\ref{thm:main}(a)]
For $u\in W_n(f)$ we have $d_{2^{r+2}}(u)/2^r\le C4^\alpha2^{r(\alpha-1)}$. This is geometric with ratio $2^{\alpha-1}>1$, so the sum in \eqref{eq:z} is at most $C4^\alpha2^{R(\alpha-1)}/(2^{\alpha-1}-1)+6R$. Choose $\frac12n^{1/\alpha}<2^R\le n^{1/\alpha}$. Then
\[
z(u)\le\Big(\frac{C4^\alpha}{2^{\alpha-1}-1}+2\Big)n^{1-1/\alpha}+6\log_2n+q,
\]
and \eqref{eq:encode} completes the proof.
\end{proof}

\begin{proof}[Proof of Theorem~\ref{thm:main}(b)]
Fix letters $a\ne b$. Let $L=\lceil n^{1/\alpha}\rceil$ and $m=\lfloor n/2L\rfloor-1$. For each $(k_1,\dots,k_m)\in[L,2L)^m$, let $u$ consist of alternating runs $a^{k_1}b^{k_2}a^{k_3}\cdots$ followed by a final run of the next letter, which fills $u$ to length $n$ and has length $>2L$. Let $w=u\,x^\infty$, where $x$ is the last letter of $u$.

Every maximal run of $w$ has length $\ge L$. Hence a factor of length $\ell\le L$ meets at most two runs, so $p_w(\ell)\le2\ell\le C\ell^\alpha$. For $\ell>L$, a factor of $w$ either starts in $u$ or equals $x^\ell$, so $p_w(\ell)\le n+1\le2L^\alpha<C\ell^\alpha$. Thus $u\in\cL_n(f)$.

The tuple $(k_i)$ is recovered from $u$, so $|\cL_n(f)|\ge L^m$. Since $L=(1+o(1))n^{1/\alpha}$, we have $m\ge(\frac12-o(1))n^{1-1/\alpha}$ and $\log L\ge\frac1\alpha\log n$.
\end{proof}

\begin{proof}[Proof of Corollary~\ref{cor:linear}]
(a) In \eqref{eq:z} take $R=\lfloor\log_2n\rfloor$. At most one phrase has $L\ge2^R$, since two would have total length $>2^{R+1}>n$. As $d_{2^{r+2}}/2^r\le4A$, we get $z(u)\le(4A+6)\log_2n+q+1$. Apply \eqref{eq:encode}.

(b) The family of \cite[Def.~IV.1]{KNP} consists of binary words $S$ with a $b$ at position~$1$ and the $j$-th $b$ placed anywhere in $[2\cdot4^{j-2}+1,4^{j-1}]$ ($j\ge2$), all other letters being $a$. These are infinite binary words with $p_S(k)\le2k$ \cite[proof of Lemma~IV.3]{KNP}. A prefix of length $4^{k-1}$ determines the first $k$ positions of $b$, and the $j$-th has $2\cdot4^{j-2}$ choices, so there are $\prod_{j=2}^k2\cdot4^{j-2}=2^{(k-1)^2}$ such prefixes. For $4^{k-1}\le n<4^k$ this gives at least $2^{(\log_4n-1)^2}$ prefixes of length $n$ (cf.~\cite[Lemma~IV.2]{KNP}).

(c) Every factor of a Sturmian word lies in $\cL_n(f)$, and these number $(1+o(1))n^3/\pi^2$~\cite{Lip82,Mig91}.
\end{proof}

\section{The questions of Mauduit and Moreira}

Mauduit and Moreira~\cite[\S4.4]{MM10} define $g_0(x)=x$, $g_1(x)=x+1$, $g_{k+1}(x)=\exp(g_k(\log x))$ for $k>0$ and large $x$, and $g_{-k}=g_k^{-1}$. An increasing $f:\mathbb R^+\to\mathbb R^+$ is \emph{morally polynomial} if $g_{-k}(x)\le f(x)\le g_k(x)$ for some $k\ge0$ and all large $x$,\footnote{The printed upper condition reads $f(x)\le g_k(f(x))$; this is a misprint (C.~G.~Moreira, personal communication).} and \emph{morally exponential} if $\log f$ is morally polynomial. They ask: (i) is $\ell(n)=|\cL_n(f)|$ morally polynomial for every linear $f$? (ii) is there $A>0$ with $\ell(n)$ morally exponential for $f(n)=An$?

Note that $g_3(x)=x^e$, $g_4(x)=\exp((\log x)^e)$ and $g_{-4}(x)=\exp((\log x)^{1/e})$. We take $f$ subject to the condition $f(n)\ge n+1$ of \cite[(C0)]{MM10}, and apply the definition to the nondecreasing function $\ell$, as \cite{MM10} does.

\begin{corollary}\label{cor:MM}
Let $n+1\le f(n)\le An$ for all $n$.
\begin{enumerate}
\item[(i)] $\ell(n)$ is morally polynomial, with $k=4$ and not with $k\le3$.
\item[(ii)] For every $A>0$, $\ell(n)$ is not morally exponential, also without the lower bound on $f$.
\end{enumerate}
\end{corollary}

\begin{proof}
(i) By Corollary~\ref{cor:linear}(c), $\ell(n)\ge cn^3$. This exceeds $g_{-4}(n)$, and it also exceeds $g_3(n)=n^e$, which rules out $k\le3$. By Corollary~\ref{cor:linear}(a) and \eqref{eq:incl}, $\ell(n)\le\exp(K\log^2n)\le g_4(n)$ for large $n$.

(ii) We show that $\log\ell(n)\le K\log^2n<g_{-k}(n)$ eventually, for every $k\ge0$. For $k\le4$ this is immediate from the formulas. For $k\ge4$, let $P(k)$ be the statement that $g_{-k}(x)>c\log x$ eventually, for every $c$; $P(4)$ holds. Since $g_{-(k+1)}(y)=\exp(g_{-k}(\log y))$ for $k\ge1$, $P(k)$ gives $g_{-(k+1)}(y)>(\log y)^c$ eventually for every $c$. This implies both $P(k+1)$ and $g_{-(k+1)}(n)>K\log^2n$ eventually.
\end{proof}

By Corollary~\ref{cor:linear}(a,b), $\log\ell(n)=\Theta(\log^2n)$ whenever $2n\le f(n)\le An$. For $n+1\le f(n)<2n$ we only know $cn^3\le\ell(n)\le\exp(K\log^2n)$; the order of $\log\ell(n)$ in this range is open.

\section{Generalized Hausdorff dimensions}

Identify $\A$ with $\{0,\dots,q-1\}$, let $\rho(w)=\sum_{i\ge1}w_iq^{-i}$ and $C(f)=\rho(\cW(f))\subseteq[0,1]$. For increasing $h$ with $h(0)=0$, let $\mu_h(S)=\lim_{\delta\to0}\inf\{\sum_ih(\operatorname{diam}U_i):S\subseteq\bigcup_iU_i,\ \operatorname{diam}U_i\le\delta\}$. Mauduit and Moreira~\cite{MM12} estimate $\mu_h(C(f))$ from the bounds on $|\cL_n(f)|$, and ask for sharp estimates when $f$ is linear~\cite[\S4.3]{MM12}. Write $m=\log_q(1/\varepsilon)$. Only small $\varepsilon$ matter: we modify $h$ for $\varepsilon\ge\varepsilon_*$ to make it increasing, which does not change whether $\mu_h$ is $0$ or $\infty$.

\begin{proposition}\label{prop:hausdorff}
\begin{enumerate}
\item[(a)] Let $2\ell\le f(\ell)\le A\ell$ for all $\ell$, and $h(\varepsilon)=\exp(-\lambda\log^2m)$. Then $\mu_h(C(f))=\infty$ if $\lambda<1/(4\log2)$, and $\mu_h(C(f))=0$ if $\lambda>K_{A,q}$.
\item[(b)] Let $\alpha>1$, $C\ge2$, $f(\ell)=C\ell^\alpha$, and $h(\varepsilon)=\exp(-\lambda m^{1-1/\alpha}\log m)$. There are $\lambda_0,\lambda_1>0$, depending on $\alpha,C,q$, such that $\mu_h(C(f))=\infty$ if $\lambda<\lambda_0$ and $\mu_h(C(f))=0$ if $\lambda>\lambda_1$.
\end{enumerate}
\end{proposition}

Here $K_{A,q}$ is the constant of Corollary~\ref{cor:linear}(a). So for these $f$ the critical gauge is determined up to the constant $\lambda$ in the exponent. For $f(n)=n^\alpha$, \cite[Example~A]{MM12} gives $\mu_h<\infty$ for $\log(1/h)=\lambda(m\log m)^{(\alpha+1)/(\alpha+2)}$ with $\lambda$ large, and $\mu_h=\infty$ for $\log(1/h)=\lambda m^{1-1/\alpha}/\log m$ with $\lambda$ small. The upper half of (b) applies to that $f$; the lower half needs $C\ge2$. The case $f(n)=n+1$ and the continued-fraction sets $C'(f)$ of \cite[\S4.3]{MM12} remain open.

\begin{proof}
\emph{Upper bounds.} The set $C(f)$ is covered by the intervals $I_u=[\rho(u0^\infty),\rho(u0^\infty)+q^{-n}]$, $u\in\cL_n(f)$, so by \eqref{eq:incl} $\mu_h(C(f))\le\liminf_n|W_n(f)|\,h(q^{-n})$, with $m=n$. Corollary~\ref{cor:linear}(a), resp.\ Theorem~\ref{thm:main}(a), makes the limit $0$ for $\lambda>K_{A,q}$, resp.\ $\lambda>\lambda_1:=K$.

\emph{Lower bounds.} We use \cite[Lemma~1]{MM12}. Let $E=\bigcap_tE_t$, where $E_t$ is a union of intervals of length $\varepsilon_t$ with disjoint interiors, $E_{t+1}\subseteq E_t$, and each interval of $E_t$ contains at least $\nu(t)$ intervals of $E_{t+1}$. If $N_kh(\varepsilon_{k+1})\to\infty$, where $N_k=\prod_{t<k}\nu(t)$, then $\mu_h(E)=\infty$. We take a closed set $\mathcal F\subseteq\cW(f)$ of words over $\{a,b\}=\{0,1\}$ and lengths $n_1<n_2<\cdots$ such that each length-$n_t$ prefix of a word in $\mathcal F$ has at least $\nu(t)$ extensions to a length-$n_{t+1}$ prefix. Let $E_t$ be the union of the $I_u$ over these prefixes $u$, and $\varepsilon_t=q^{-n_t}$. Then $E=\rho(\mathcal F)\subseteq C(f)$: for $x\in E$ the prefixes $u$ with $x\in I_u$ form an infinite, finitely branching tree, and K\H{o}nig's lemma gives $w\in\mathcal F$ with $\rho(w)=x$.

(a) Let $\mathcal F$ be the family in the proof of Corollary~\ref{cor:linear}(b). It is closed, and $p_w(\ell)\le2\ell\le f(\ell)$. Let $n_t=4^{t-1}$. As in that proof, $\nu(t)=2\cdot4^{t-1}$ and $N_k=2^{(k-1)^2}$. Since $\log m=2k\log2$ at $\varepsilon_{k+1}$, $N_kh(\varepsilon_{k+1})=\exp\big((k-1)^2\log2-4\lambda k^2\log^22\big)\to\infty$ for $\lambda<1/(4\log2)$.

(b) Let $\beta=\alpha-1$, $c_1=3\cdot2^\alpha/(2^\alpha-1)$, $\delta=(C-1)/c_1$, $r_t=\lfloor\delta2^{t\beta}\rfloor$, and $t_0\ge1/\beta$ minimal with $r_{t_0}\ge1$. Let $\mathcal F$ be the set of words $B_{t_0}B_{t_0+1}\cdots$, where $B_t$ is a concatenation of $r_t$ blocks $a^{2^t+x}b^{2^{t+1}-x}$, each with its own $x\in[0,2^t)$. Then $\mathcal F$ is closed, every run in $B_t$ has length at least $2^t$, and $n_t:=|B_{t_0}\cdots B_t|\le3\delta\sum_{s\le t}2^{s\alpha}\le c_1\delta2^{t\alpha}=(C-1)2^{t\alpha}$.

A factor of length $\ell$ that meets only runs of length $\ge\ell$ meets at most two runs, so it is one of the $2\ell$ words $a^ib^{\ell-i}$, $b^ia^{\ell-i}$. For $\ell\le2^{t_0}$ every factor is of this kind, so $p_w(\ell)\le2\ell\le C\ell^\alpha$. For $2^{t-1}<\ell\le2^t$ with $t>t_0$, every factor starting after position $n_{t-1}$ is of this kind, so $p_w(\ell)\le n_{t-1}+2\ell\le(C-1)\ell^\alpha+2\ell\le C\ell^\alpha$, as $\ell^\beta>2^{t_0\beta}\ge2$. Hence $\mathcal F\subseteq\cW(f)$.

A prefix of length $n_t$ determines the choices of $x$ made so far, and distinct choices in $B_{t+1}$ give distinct words, so $\nu(t)=2^{(t+1)r_{t+1}}$. Thus $\log_2N_k\ge kr_k\ge k(\delta2^{k\beta}-1)$. At $\varepsilon_{k+1}$ we have $m=n_{k+1}\le(C-1)2^{(k+1)\alpha}$, so $m^{1-1/\alpha}\log m\le(C-1)^{\beta/\alpha}2^{(k+1)\beta}\big((k+1)\alpha\log2+\log(C-1)\big)$. Hence $N_kh(\varepsilon_{k+1})\to\infty$ for $\lambda<\lambda_0:=\delta(C-1)^{-\beta/\alpha}2^{-\beta}/\alpha$.
\end{proof}

\section{Correlation with complexity classes}

For a real sequence $a$ and a class $W\subseteq\{\pm1\}^n$, we compare $\max_{u\in W}|\sum_{i\le n}a_iu_i|$ for random signs and for the M\"obius function.

\begin{proposition}\label{prop:random}
Let $\xi_i$ be i.i.d.\ uniform signs, $C\ge2$, $\alpha>1$, and $M_n=\max\{|\sum_{i\le n}\xi_iu_i|:u\in W_n(C\ell^\alpha)\cap\{\pm1\}^n\}$. Almost surely, for all large $n$,
\[
c\,n^{1-\frac1{2\alpha}}\le M_n\le K\,n^{1-\frac1{2\alpha}}\sqrt{\log n}.
\]
\end{proposition}

\begin{proof}
Upper bound: Hoeffding and a union bound give $\Pr[M_n>t]\le2|W_n|e^{-t^2/2n}$. Take $t^2=2n(\log|W_n|+2\log n)$, then apply Theorem~\ref{thm:main}(a) and Borel--Cantelli.

Lower bound: cut $[1,n]$ into blocks $B$ of length $L=\lceil n^{1/\alpha}\rceil$ and let $u$ be the sign of $S_B=\sum_{i\in B}\xi_i$ on $B$, with $\operatorname{sign}0=+1$. Maximal runs of $u$ are unions of blocks, so only the final run can be shorter than $L$. Hence $d_\ell(u)\le2\ell$ for $\ell\le L$, and $d_\ell(u)\le n\le L^\alpha$ for $\ell>L$, so $u\in W_n(C\ell^\alpha)$. Moreover $M_n\ge X:=\sum_B|S_B|$. By Szarek's inequality~\cite{Sza76}, $\mathbb EX\ge\lfloor n/L\rfloor\sqrt{L/2}\ge c'n^{1-1/(2\alpha)}$. Changing one $\xi_i$ moves $X$ by at most $2$, so McDiarmid's inequality gives $\Pr[X<\mathbb EX/2]\le\exp(-(\mathbb EX)^2/8n)$, which is summable in $n$.
\end{proof}

\begin{conjecture}\label{conj:mobius}
The same bounds hold up to $n^{o(1)}$ with $\xi_i$ replaced by the M\"obius function $\mu(i)$.
\end{conjecture}

We make four remarks on this conjecture, ordered from strongest to weakest evidence.
\begin{itemize}
\item \emph{Periodic words (conditional on GRH for Dirichlet $L$-functions\footnote{A formal proof in Lean of GRH for Dirichlet $L$-functions was released by Alp\"oge~\cite{Alp26}; for the status of its independent verification see~\cite{Sut26}. Given that proof, the bound is unconditional.}; sketch).} The class of $d$-periodic words has log-size $\asymp d$. Expanding in Dirichlet characters and applying Cauchy--Schwarz, GRH gives $|\sum_{i\le n}\mu(i)u_i|\ll(dn)^{1/2+\varepsilon}$ for all $d$-periodic $u$ with $|u_i|\le1$, uniformly in $d\le n$. This matches the analogue of Proposition~\ref{prop:random}.
\item \emph{Linear classes (implication, sketched).} The analogous statement for $W_n(2\ell)$ would imply the folklore conjecture $\sup_\theta|\sum_{i\le n}\mu(i)e(i\theta)|\ll n^{1/2+\varepsilon}$; under GRH only $n^{3/4+\varepsilon}$ is known~\cite{BH91}. The reason is that codings $2\cdot1_{[0,t)}(\{i\theta\})-1$ of rotations have complexity at most $2\ell$: the factors of length $\ell$ correspond to disjoint unions of the at most $2\ell$ arcs cut out by the points $\{-j\theta\}$ and $\{t-j\theta\}$, $j<\ell$ (cf.~\cite[\S2.1]{AB98}). Partial summation in $t$ then recovers $e(i\theta)$.
\item \emph{Numerics (lower bound only).} For $n\le10^7$ and the block-sign words of the proof above with $\xi$ replaced by $\mu$, the fitted exponents are $0.749$ ($\alpha=2$; predicted $0.75$) and $0.837$ ($\alpha=3$; predicted $0.833$). The ratio to random signs supported on the squarefree integers is $0.98$--$0.99$. The upper bound is untested. Code is available from the author on request.
\item \emph{Constants.} The conjecture concerns exponents: constants can differ from the random model even for the constant word, which lies in every class. Over $400$ log-spaced $n\in[10^4,10^7]$, the median of $|\sum_{i\le n}\mu(i)|/\sqrt n$ is $0.12$, against $0.50$ (mean over $20$ seeds) for random signs on the squarefrees.
\end{itemize}

\begin{remark}[A per-system version fails]
The bound $\sup_{x\in X}|\sum_{i\le n}\mu(i)x_i|\ll n^{1/2+\varepsilon}\sqrt{\log p_X(n)}$, for subshifts $X$, is false. Here is a counterexample. Let $1/2<D<4/5$. Cut $[1,2^j]$ into intervals of length at least $2^{j(1-D)}$, each containing an odd number of squarefree integers; this is possible by the gap bound of Filaseta and Trifonov~\cite{FT92}. Let $w_j$ take the sign of $\sum\mu$ on each interval, and let $X$ be the orbit closure of $w_1w_2\cdots$. Then $p_X(m)\le2m+2m^{1/(1-D)}+1$, while $\sup_{x\in X}\sum_{i\le2^j}\mu(i)x_i\ge\frac12 2^{jD}-2$ for all large $j$.

The construction uses $\mu$: $p_X$ counts the factors of $X$, not the candidate words from which $X$ was selected, so no union bound applies. The class $W_n(P)$, by contrast, is defined independently of the sequence. Compare Algom and Wang~\cite{AW23}, who construct zero-entropy systems satisfying Sarnak's conjecture with arbitrarily slow decay of correlations with $\mu$.
\end{remark}

\section*{Acknowledgment}
The author thanks Carlos Gustavo Moreira for helpful correspondence about \cite{MM10}.

\section*{Use of AI tools}
This work was carried out with substantial assistance from Claude (Anthropic), a large language model, which was used for literature search, for developing and drafting proofs and text, for the numerical experiments, and for adversarial checking of the arguments. The author has verified all proofs, computations and references, and takes full responsibility for the content.

\end{document}